\documentclass[12pt,francais,nothms]{aclart}
\usepackage{mathptm}
\usepackage[all]{xy}

\newtheorem{thm}{Th\'eor\`eme}

\usepackage[all]{xy}

\RequirePackage[T1]{fontenc}
\RequirePackage{mathrsfs,amsmath,amssymb}
\let\mathcal\mathscr
 
\RequirePackage{url}

\RequirePackage{xspace}
\RequirePackage[english,frenchb]{babel}

\theoremstyle{remark}

\def\ord{{\rm ord}}
\def\ordjac{{\rm ordjac}}

\def\11{{\mathbf 1}}

\def\Spec{\operatorname{Spec}}

\def\Def{{\rm Def}}
\def\GDef{{\rm GDef}}
\def\RDef{{\rm RDef}}

\def\LPas{\cL_{\rm DP}}

\def\AA{{\mathbf A}}

\def\LL{{\mathbf L}}

\def\NN{{\mathbf N}}

\def\QQ{{\mathbf Q}}

\def\ZZ{{\mathbf Z}}

\def\cA{{\mathcal A}}

\def\cC{{\mathcal C}}

\def\cL{{\mathcal L}}
\def\cM{{\mathcal M}}

\def\cO{{\mathcal O}}

\def\cU{{\mathcal U}}

\def\cX{{\mathcal X}}
\def\cY{{\mathcal Y}}

\mathchardef\alphag="7C0B
\mathchardef\betag="7C0C
\mathchardef\gammag="7C0D
\mathchardef\deltag="7C0E
\mathchardef\varepsilong="7C22
\mathchardef\varphig="7C27
\mathchardef\psig="7C20
\mathchardef\zetag="7C10
\mathchardef\epsilong="7C0F
\mathchardef\rhog="7C1A
\mathchardef\taug="7C1C
\mathchardef\upsilong="7C1D
\mathchardef\iotag="7C13
\mathchardef\thetag="7C12
\mathchardef\pig="7C19
\mathchardef\sigmag="7C1B
\mathchardef\etag="7C11
\mathchardef\omegag="7C21
\mathchardef\kappag="7C14
\mathchardef\lambdag="7C15
\mathchardef\mug="7C16
\mathchardef\xig="7C18
\mathchardef\chig="7C1F
\mathchardef\nug="7C17
\mathchardef\varthetag="7C23
\mathchardef\varpig="7C24
\mathchardef\varrhog="7C25
\mathchardef\varsigmag="7C26
\mathchardef\Omegag="7C0A
\mathchardef\Thetag="7C02
\mathchardef\Sigmag="7C06
\mathchardef\Deltag="7C01
\mathchardef\Phig="7C08
\mathchardef\Gammag="7C00
\mathchardef\Psig="7C09
\mathchardef\Lambdag="7C03
\mathchardef\Xig="7C04
\mathchardef\Pig="7C05
\mathchardef\Upsilong="7C07

\begin{document}

\selectlanguage{french}
\title[Fonctions constructibles et int\'egration motivique II]
  {{\normalsize\bfseries G\'eom\'etrie alg\'ebrique/\itshape Algebraic geometry \\}
  
Fonctions constructibles et int\'egration motivique II}
\alttitle{Constructible functions and motivic integration II}

\author{Raf Cluckers}
\address{Katholieke Universiteit Leuven, Department of Mathematics,
Celestijnenlaan 200B, 3001 Leu\-ven, Bel\-gium }
\email{raf.cluckers@wis.kuleuven.ac.be}

\author{Fran\c cois Loeser}

\address{{\'E}cole Normale Sup{\'e}rieure,
D{\'e}partement de math{\'e}matiques et applications,
45 rue d'Ulm,
75230 Paris Cedex 05, France
(UMR 8553 du CNRS)}
\email{Francois.Loeser@ens.fr}
\date{\today}

\begin{altabstract}We extend the formalism of \cite{cr1}
to a global setting 
for which 
a theorem on fiber integrals and a Fubini
theorem are obtained.
We compare our formalism to the previous constructions 
given in \cite{arcs} and  \cite{JAMS}.
Details of constructions and proofs will
be given in \cite{cl}.
\end{altabstract}

\begin{abstract}
On \'etend le formalisme de \cite{cr1} \`a un cadre global
dans lequel on \'etablit 
un th\'eor\`eme d'int\'egration dans les fibres
ainsi qu'un th\'eor\`eme de Fubini.
On compare notre formalisme 
avec les constructions ant\'erieures de \cite{arcs} et de \cite{JAMS}.
Les d\'etails des constructions et des preuves seront donn\'es
dans \cite{cl}.

\end{abstract}

\maketitle

\selectlanguage{french}

\section{Formes volume sur les sous-assignements d\'efinissables globaux}
\subsection{Sous-assignements d\'efinissables globaux}Dans
ce travail on se place dans le  cadre de  \cite{cr1} dont on reprend les
notations.
En particulier
on se donne un corps
$k$ de caract\'eristique $0$ et on note
$F_k$ la cat\'egorie des corps contenant $k$.
Dans la note  \cite{cr1} on a d\'efini la  cat\'egorie
$\Def_k$
des sous-assignements d\'efinissables.
Les objets de $\Def_k$ sont de nature affine, 
\'etant des sous-assignements
des foncteurs $h [m, n, r] : F_k \rightarrow {\rm Ens}$ donn\'es
par $K \mapsto K((t))^m \times K^n \times \ZZ^r$. 
On va maintenant consid\'erer leur analogue global.

Soit $\cX$ une vari\'et\'e, c'est \`a dire un sch\'ema de type fini s\'epar\'e et r\'eduit,
sur $k ((t))$ et soit $X$ une vari\'et\'e sur $k$. Pour $r$ un entier $\geq 0$
on note $h [\cX, X, r]$ le foncteur 
$F_k \rightarrow {\rm Ens}$ donn\'e
par $K \mapsto \cX (K ((t))) \times X (K) \times \ZZ^r$.
Lorsque $X = \Spec k$ et $r = 0$, on \'ecrit
$h [\cX]$ pour $h [\cX, X, r]$.
Si $\cX$ et $X$ sont affines et si $i : \cX \hookrightarrow \AA^m_{k ((t))}$
et $j : X \hookrightarrow \AA^n_k$ sont des immersions ferm\'ees,
on dit qu'un sous-assignement $h$ de $h [\cX, X, r]$ est d\'efinissable
si son image par le morphisme
$h [\cX, X, r] \rightarrow h [m, n, r]$ associ\'e \`a $i$ et $j$ est un sous-assignement
d\'efinissable de $h [m, n, r]$. Cette d\'efinition ne d\'epend pas du choix de $i$
et de $j$. En g\'en\'eral, on dira qu'un sous-assignement $h$
de $h [\cX, X, r]$ est d\'efinissable s'il existe des recouvrements
$(\cU_i)$ et $(U_j)$ de $\cX$ et $X$ par des ouverts affines tels 
$h \cap h [\cU_i, U_j, r]$ soit un sous-assignement
d\'efinissable de $h [\cU_i, U_j, r]$ pour tout $i$ et $j$.
On obtient ainsi une cat\'egorie
$\GDef_k$ dont les objets
sont les 
sous-assignements
d\'efinissables d'un $h [\cX, X, r]$, les morphismes
\'etant les morphismes d\'efinissables, c'est \`a dire ceux dont le graphe
est un sous-assignement
d\'efinissable.

La cat\'egorie $\Def_k$ est une sous-cat\'egorie pleine
de $\GDef_k$. La notion de dimension introduite dans \cite{cr1} s'\'etend imm\'ediatement
aux objets de $\GDef_k$, ainsi que le fait que la dimension
soit invariante par isomorphisme (Proposition 1 de \cite{cr1}).
Si $S$ est un
objet de $\GDef_k$, les d\'efinitions de la cat\'egorie
$\RDef_S$, du semi-anneau
$\cC_+ (S)$, de
l'anneau $\cC (S)$,
du $\cC_+ (S)$-semi-module gradu\'e $C_+ (S)$ et du
$\cC (S)$-module gradu\'e $C (S)$ donn\'ees en 
\cite{cr1} s'\'etendent directement.

\subsection{Formes diff\'erentielles d\'efinissables}
Soit $h$ un sous-assignement
d\'efinissable d'un $h [\cX, X, r] $. On note
$\cA (h)$ l'anneau des morphismes d\'efinissables
$h  \rightarrow
h [\AA^1_{k ((t))}]$. On va d\'efinir, pour $i$ dans
$\NN$, le  $\cA (h)$-module
$\Omega^i (h)$ des $i$-formes d\'efinissables sur $h$.
Soit $\cY$ le  ferm\'e  de $\cX$, adh\'erence de Zariski de l'image de  $h$  
par la  projection $\pi : h  [\cX, X, r] \rightarrow
h [\cX] $. On consid\`ere
le faisceau $\Omega^i_{\cY}$ des $i$-formes alg\'ebriques sur
$\cY$, on note $\cA_{\cY}$ le faisceau Zariski
$U \mapsto \cA (h [U])$
sur $\cY$, $\Omega^i_{h [\cY]}$
le faisceau 
$\cA_{\cY} \otimes_{\cO_{\cY}} \Omega^i_{\cY}$,
et on pose
$\Omega^i (h) :=
\cA (h) \otimes_{\cA (h [\cY])}    \Omega^i_{h [\cY]} (\cY)$, la 
structure de $\cA (h [\cY])$-alg\`ebre sur 
$\cA (h)$ \'etant donn\'ee par  composition avec $\pi$.

\subsection{Formes volume d\'efinissables}On suppose maintenant que $h$ est de dimension $d$.
On note $\cA^< (h)$ l'id\'eal des fonctions dans $\cA (h)$
qui sont nulles en dehors d'un sous-assignement d\'efinissable de
dimension $<d$.
On a un morphisme canonique de groupes ab\'eliens
$\lambda : \cA (h) / \cA^< (h) \rightarrow C^d_+ (h)$
qui envoie la classe d'une fonction $f$ sur  
la classe de $\LL^{- \ord f}$,
avec la convention
$\LL^{- \ord 0} = 0$.
On pose
$\tilde \Omega^d (h) = \cA (h) / \cA^< (h) \otimes_{\cA^< (h)} \Omega^d (h) $,
et on d\'efinit l'ensemble $\vert \Omega \vert_+ (h) $
des formes volume positives d\'efinissables
comme 
le quotient du semi-groupe ab\'elien libre sur les symboles $(\omega, g)$
avec $\omega$ dans $\tilde \Omega^d (h) $ et $g$ dans $C^d_+ (h)$
par les relations
$(f \omega, g) = (\omega, \lambda (f) g)$, $(\omega, g + g') = 
(\omega, g) + (\omega, g')$ et $(\omega, 0)= 0$.
On \'ecrira $g |\omega|$ pour la classe de
$(\omega, g)$, de fa\c con \`a ce que
$g \vert f \omega \vert = g \LL^{- \ord f} \vert \omega \vert$.
La structure de $\cC_+ (h)$-semi-module sur
$C^d_+ (h)$ induit par passage au quotient une structure
de semi-anneau 
sur $C^d_+ (h)$ et
$\vert \tilde \Omega (h) \vert_+$
est muni naturellement d'une structure de
$C^d_+ (h)$-semi-module.
On dit que $\vert \omega \vert$ dans $\vert \tilde \Omega (h) \vert_+$
est une forme jauge si c'est un g\'en\'erateur de ce 
semi-module. 
On v\'erife qu'il existe toujours des formes jauges.
On d\'efinit de fa\c con similaire 
$\vert \tilde \Omega \vert(h)$
en rempla\c cant $C^{d}_+$ par $C^{d}$.
Pour des raisons de place on ne donnera toutefois dans cette note
que des \'enonc\'es concernant les formes volumes
positives.

Si $h$ est un sous-assignement d\'efinissable de dimension $d$ de
$h [m, n, r]$, on dispose, de fa\c con analogue \`a la construction de Serre
\cite{serre} dans le cas  $p$-adique, d'une forme 
jauge canonique sur $h$, not\'ee $\vert \omega_0\vert_h$.
D\'esignons par $x_1$, \dots, $x_m$ les
coordonn\'ees sur $\AA^m_{k((t))}$ et consid\'erons les  $d$-formes
$\omega_I := dx_{i_1} \wedge \dots \wedge dx_{i_d}$ pour $I =
\{i_1, \dots, i_d\} \subset \{1, \dots m\}$, $i_1 < \dots < i_d$
et leur image $\vert \omega_I \vert_h$ dans
$\vert \Omega \vert_+ (h)$. On v\'erifie qu'il existe un unique 
\'el\'ement  $\vert \omega_0\vert_h$ de 
$\vert \tilde \Omega \vert_+ (h)$, tel que, pour tout $I$, il existe
une fonction d\'efinissable
\`a valeurs enti\`eres $\alpha_I$ sur 
$h$ telle que $\vert \omega_I \vert_h =
\LL^{- \alpha_I} \vert \omega_0\vert_h$
dans $\vert \tilde \Omega \vert_+ (h)$ et telle que 
$\inf_I \alpha_I = 0$ en dehors d'un sous-assignement
d\'efinissable de dimension $<d$.

Si $f : h \rightarrow h'$ est un morphisme dans $\GDef_k$
avec $h$ et $h'$ de dimension $d$
et dont toutes les fibres sont de dimension $0$, on dispose
d'une application $f^* : \vert \tilde \Omega \vert_+ (h')  \rightarrow \vert \tilde \Omega \vert_+ (h)$
induite par le 
pull-back des formes diff\'erentielles.
Ceci r\'esulte du fait que 
$f$ est ``analytique''
dans le compl\'ementaire d'un sous-assignement d\'efinissable de
dimension $d - 1$ de $h$.
Si $h$ et $h'$ sont de plus des objets de $\Def_k$,
on a $f^* \vert \omega_0\vert_{h'} =
\LL^{ -\ordjac f} \vert \omega_0\vert_h$, avec
$\ordjac f$ l'ordre du jacobien de $f$ (cf. \cite{cr1}).

Si $\cX$ est une $k((t))$-vari\'et\'e de dimension $d$,
admettant un $k [[t]]$-mod\`ele $\cX^0$,
on peut d\'efinir un \'el\'ement $|\omega_0|$ de
$|\tilde \Omega|_+ (h [\cX])$, ne d\'ependant que de
$\cX^0$, et  caract\'eris\'e par le fait que pour tout ouvert
$U^0$ de $\cX^0$ sur lequel 
le $k [[t]]$-module
$\Omega^d_{U^0 | k [[t]]} (U^0)$ est engendr\'e
par une forme $\omega$, $|\omega_0|_{|h [U^0 \otimes \Spec k((t))]} = |\omega|$
dans $|\tilde \Omega|_+(h [U^0 \otimes \Spec k((t))])$. 
 
\section{Int\'egration des formes volume et th\'eor\`eme de Fubini}
Soit $f : S \rightarrow S'$ un morphisme dans $\Def_k$,
avec $S$ de dimension $s$ et $S'$ de dimension $s'$.
Toute forme volume positive
$\alpha$ dans
$|\tilde \Omega|_+ (S)$ s'\'ecrit $\alpha = \psi_{\alpha} |\omega_0|_S$
avec $\psi_{\alpha}$ dans $C^s_+ (S)$. 
On dit que $\alpha$ est $f$-int\'egrable
si 
$\psi_{\alpha}$
est $f$-int\'egrable
et on pose alors
$$f_!^{\rm top} (\alpha) := \{f_! (\psi_{\alpha})\}_{s'}|\omega_0|_{S'},$$
$\{f_! (\psi_{\alpha})\}_{s'}$ d\'esignant la composante de
$f_! (\psi_{\alpha})$ dans $C^{s'}_+ (S')$.
Consid\'erons maintenant 
$f : S \rightarrow S'$ un morphisme dans $\GDef_k$.
Supposons qu'il existe des isomorphismes
$\varphi : T \rightarrow S$ et $\varphi' : T' \rightarrow S'$
avec $T$ et $T'$ dans $\Def_k$.
On note $\tilde f$ le morphism $ T \rightarrow T'$
v\'erifiant $\varphi' \circ \tilde f = f \circ \varphi$.
On dira que $\alpha$ dans
$|\tilde \Omega|_+ (S)$ est $f$-int\'egrable si 
$\varphi^* (\alpha)$ est $\tilde f$-int\'egrable
et on d\'efinit alors $f^{\rm top}_! (\alpha)$
par la relation
$$
{\tilde f}^{\rm top}_! (\varphi^* (\alpha))
= \varphi'{}^* (f^{\rm top}_! (\alpha)).
$$
Il r\'esulte
du th\'eor\`eme 2 de \cite{cr1} que
cette d\'efinition ne d\'epend
pas du choix des isomorphismes $\varphi$ et $\varphi'$.
Par additivit\'e, et en utilisant
des cartes affines, on \'etend ce qui pr\'ec\`ede
\`a un morphisme arbitraire
$f : S \rightarrow S'$ dans 
$\GDef_k$, de fa\c con \`a d\'efinir 
les formes volumes $f$-int\'egrables dans 
$|\tilde \Omega|_+ (S)$  et pour de telles formes
$\alpha$ l'int\'egrale dans les fibres
$f^{\rm top}_! (\alpha)$ appartenant \`a 
$|\tilde \Omega|_+ (S')$.
Quand $S = h [0, 0, 0]$, on dit que $\alpha$ est
int\'egrable au lieu de $f$-int\'egrable, et on \'ecrit
$\int_S \alpha$ \`a la place de
$f^{\rm top}_! (\alpha)$.


Le r\'esultat suivant est une cons\'equence du th\'eor\`eme 1 de
\cite{cr1}.

\begin{thm}[Th\'eor\`eme de Fubini]\label{if}Soit $f : S \rightarrow S'$
un morphisme dans $\GDef_k$.
On suppose que $S$ est de dimension $s$, que $S'$ est de dimension $s'$,
et que toutes les fibres de $f$ sont de dimension $s - s'$.
Une forme volume positive $\alpha$ dans $|\tilde \Omega|_+ (S)$
est int\'egrable si et seulement si elle est $f$-int\'egrable et
$f^{\rm top}_! (\alpha)$ est int\'egrable.
Si c'est le cas, alors
$$
\int_S \alpha = \int_{S'} f^{\rm top}_! (\alpha).
$$
\end{thm}

\section{Comparaison avec les constructions ant\'erieures}

\subsection{Comparaison avec la construction classique}\label{3.1}Dans la d\'efinition de $\Def_k$, $\RDef_k$ et
$\GDef_k$, au lieu de consid\'erer la cat\'egorie
$F_k$ des corps contenant $k$, on pourrait choisir de se restreindre
\`a la sous-cat\'egorie ${\rm ACF}_k$ des corps
alg\'ebriquement clos contenant $k$ et d\'efinir ainsi des
cat\'egories
$\Def_{k, {\rm ACF}_k}$, etc.
Il r\'esulte du th\'eor\`eme de constructibilit\'e de Chevalley que
$K_0 (\RDef_{k, {\rm ACF}_k})$
n'est autre que l'anneau de Grothendieck
$K_0 ({\rm Var}_k)$ 
consid\'er\'e dans  \cite{arcs}.
On dispose ainsi d'un morphisme canonique
$SK_0 (\RDef_{k}) \rightarrow K_0 ({\rm Var}_k)$
envoyant $\LL$ sur la classe de $\AA^1_k$ \'egalement not\'ee
$\LL$,
que l'on peut \'etendre en un morphisme
$\gamma : SK_0 (\RDef_k) \otimes_{\NN [\LL - 1]} A_+
\rightarrow K_0 ({\rm Var}_k) \otimes_{\ZZ [\LL]} A$.
En consid\'erant le d\'eveloppement en s\'erie de
$(1 - \LL^{-i})^{-1}$, on d\'efinit
\'egalement un morphisme
canonique $\delta : K_0 ({\rm Var}_k)\otimes_{\ZZ [\LL]} A \rightarrow 
\widehat \cM$, avec $\widehat \cM$ le compl\'et\'e 
de $K_0 ({\rm Var}_k) [\LL^{-1}]$
consid\'er\'e dans \cite{arcs}.

Soit $X$ une vari\'et\'e alg\'ebrique
sur $k$ de dimension $d$.
Posons $\cX^0 := X \otimes_{\Spec k} \Spec k [[t]]$
et
$\cX := \cX^0 \otimes_{\Spec k [[t]]} \Spec k ((t))$.
Consid\'erons un sous-assignement d\'efinissable
$W$ de  $h [\cX]$ dans le langage $\LPas$, avec la restriction
que les constantes
dans la sorte de type ${\rm Val}$
apparaissant dans les 
formules d\'efinissant $W$ dans des cartes affines
d\'efinies sur $k$
sont dans $k$ (et non dans $k ((t))$).
On suppose que $W (K) \subset \cX (K [[t]])$
pour tout $K$ dans $F_k$.
Avec les notations de 
\cite{arcs}, les formules d\'efinissant 
$W$ dans une carte affine
d\'efinissent 
un sous-ensemble semi-alg\'ebrique de
l'espace des arcs
$\cL (X)$ dans la carte correspondante, 
et de cette fa\c con on associe canoniquement 
\`a
$W$ une partie semi-alg\'ebrique $\tilde W$ de
$\cL (X)$. De m\^eme \`a toute 
fonction d\'efinissable \`a valeurs enti\`eres
$\alpha$ sur $W$ (v\'erifiant la condition
additionnelle que les constantes
dans la sorte de type ${\rm Val}$
apparaissant dans les formules d\'efinissant $\alpha$
ne peuvent \^etre prises que dans $k$)
on associe une 
fonction semi-alg\'ebrique
$\tilde \alpha$ sur $\tilde W$.

\begin{thm}\label{compvar}Sous les hypoth\`eses pr\'ec\'edentes,
si $|\omega_0|$ est la forme volume canonique sur $h [\cX]$,
pour toute 
fonction d\'efinissable \`a valeurs enti\`eres
$\alpha$ sur $W$ v\'erifiant les conditions pr\'ec\'edentes
et born\'ee inf\'erieurement,
$\11_W \LL^{-\alpha}|\omega_0|$ est int\'egrable sur
$h [\cX]$
et
$$
(\delta \circ \gamma) \Bigl(\int_{h [\cX]} \11_W \LL^{-\alpha}|\omega_0|\Bigr)
=
\int_{\tilde W} \LL^{- \tilde \alpha} d \mu',
$$
$\mu'$ d\'esignant la mesure motivique construite dans  \cite{arcs}.
\end{thm}

Il est int\'eressant de remarquer qu'il r\'esulte du th\'eor\`eme \ref{compvar}
que pour
les ensembles et les fonctions semi-alg\'ebriques
l'int\'egrale motivique
construite dans \cite{arcs} existe d\'ej\`a dans
$K_0 ({\rm Var}_k) \otimes_{\ZZ [\LL]} A$, sans qu'il soit n\'ecessaire
de compl\'eter d'avantage l'anneau de Grothendieck.

\subsection{Comparaison avec l'int\'egration motivique arithm\'etique}
A la place de ${\rm ACF}_k$, on peut aussi consid\'erer la
cat\'egorie ${\rm PFF}_k$
des corps pseudo-finis contenant $k$. Rappelons qu'un tel corps
est un corps parfait $F$
admettant une unique extension de degr\'e $n $ pour tout $n$ dans une
cl\^oture alg\'ebrique fix\'ee, et tel que toute vari\'et\'e
g\'eom\'etriquement irr\'eductible sur $F$ admette un point $F$-rationnel.
Par restriction de $F_k$ \`a
${\rm PFF}_k$ on d\'efinit des
cat\'egories
$\Def_{k, {\rm PFF}_k}$, etc.
En particulier l'anneau de Grothendieck
$K_0 (\RDef_{k, {\rm PFF}_k})$
est identique \`a l'anneau not\'e
$K_0 ({\rm PFF}_k)$ dans  \cite{pek} et \cite{dw}.

Dans l'article \cite{JAMS}, 
l'int\'egrale motivique arithm\'etique prenait ses
valeurs dans 
un certain compl\'et\'e
$\hat K_0^v ({\rm Mot}_{k, \bar \QQ})_{\QQ}$ d'un anneau not\'e
$K_0^v ({\rm Mot}_{k, \bar \QQ})_{\QQ}$. 
Ce n'est qu'un peu plus tard qu'il a \'et\'e remarqu\'e
dans
\cite{pek} et \cite{dw} que l'on peut se restreindre \`a un anneau 
plus petit not\'e
$K_0^{\rm mot} ({\rm Var}_k) \otimes \QQ$,
dont nous allons maintenant rappeler la d\'efinition.

Le corps 
$k$ \'etant de caract\'eristique z\'ero, il existe, d'apr\`es 
\cite{GS}
et \cite{GN}, 
un unique morphisme d'anneaux
$K_0 ({\rm Var}_k) \rightarrow K_0 ({\rm CHMot}_k)$
associant \`a la classe
d'une vari\'et\'e $X$ projective et lisse sur $k$
la classe de son motif de Chow.
Ici
$K_0 ({\rm CHMot}_k)$
d\'esigne l'anneau de Grothendieck 
de la cat\'egorie des motifs de Chow sur $k$
(avec coefficients rationnels).
Par d\'efinition
$K_0^{\rm mot} ({\rm Var}_k) $ est
l'image de 
$K_0 ({\rm Var}_k)$ dans 
$K_0 ({\rm CHMot}_k)$ par ce morphisme.
[Noter que la d\'efinition de $K_0^{\rm mot} ({\rm Var}_k) $ 
donn\'ee dans  \cite{pek} est incorrecte.]
Dans 
\cite{pek} et \cite{dw}, les auteurs ont construit, en utilisant des 
r\'esultats de \cite{JAMS},
un morphisme canonique
$\chi_c : K_0
({\rm PFF}_k ) \rightarrow K_0^{\rm mot} ({\rm Var}_{k}) \otimes \QQ$.

La mesure motivique arithm\'etique 
prend ses valeurs dans une certaine compl\'etion
$ \hat K_0^{\rm mot} ({\rm Var}_{k}) \otimes \QQ$
de la localisation de 
$K_0^{\rm mot} ({\rm Var}_{k}) \otimes \QQ$
par rapport \`a la classe de l'image de la droite affine. 
On dispose d'un morphisme canonique
$\tilde \gamma :
SK_0 (\RDef_{k}) \otimes_{\NN [\LL- 1]} A_+
\rightarrow
K_0 ({\rm PFF}_k) \otimes_{\ZZ [\LL]} A$,
En composant avec 
$\chi_c$ et en d\'eveloppant en s\'erie $(1 - \LL^{-i})^{-1}$,
on obtient un  
morphisme canonique
$\tilde \delta :
K_0 ({\rm PFF}_k) \otimes_{\ZZ [\LL]} A
\rightarrow \hat K_0^{\rm mot} ({\rm Var}_{k} )\otimes \QQ$.

Soit $X$ une vari\'et\'e alg\'ebrique
sur $k$ de dimension $d$.
Posons $\cX^0 := X \otimes_{\Spec k} \Spec k [[t]]$,
$\cX := \cX^0 \otimes_{\Spec k [[t]]} \Spec k ((t))$
et
consid\'erons un sous-assignement d\'efinissable
$W$ de  $h [\cX]$
v\'erifiant les conditions de \ref{3.1}.
Les formules d\'efinissant 
$W$ dans une carte affine
permettent de d\'efinir, dans la terminologie et avec les notations 
de 
\cite{JAMS},
un sous-assignement 
d\'efinissable de $h_{\cL (X)}$
dans la carte correspondante, ce qui permet
d'associer canoniquement 
\`a
$W$ un sous-assignement 
d\'efinissable $\tilde W$ de $h_{\cL (X)}$.
\begin{thm}\label{arcompvar}Sous les  hypoth\`eses et avec les notations 
pr\'ec\'edentes,
$\11_W |\omega_0|$ est  int\'egrable sur 
$h [\cX]$ et
$$
(\tilde \delta \circ \tilde \gamma) \Bigl(\int_{h [\cX]} \11_W |\omega_0|\Bigr)
=
\nu (\tilde W),
$$
$\nu$ d\'esignant la mesure motivique arithm\'etique d\'efinie dans
\cite{JAMS}.
\end{thm}

Il r\'esulte du th\'eor\`eme \ref{arcompvar}
que dans le pr\'esent contexte
l'int\'egrale motivique arithm\'etique
construite dans \cite{JAMS} existe d\'ej\`a dans
$K_0 ({\rm PFF}_k) \otimes_{\ZZ [\LL]} A$, sans qu'il soit n\'ecessaire
de compl\'eter d'avantage l'anneau de Grothendieck ni de passer
aux motifs de Chow.

\subsection*{}{\small Pendant la r\'ealisation de ce projet, 
le premier auteur  \'etait chercheur postdoctoral du Fonds 
de Recherche Scientifique - Flandres (Belgique)
et il a b\'en\'efici\'e du soutien partiel du projet europ\'een
EAGER.}

 \bibliographystyle{smfplain}
 \bibliography{aclab,acl}
\end{document}